\documentclass[11pt,
]{article}

\usepackage[tbtags]{amsmath}
\usepackage{amsfonts,amssymb,mathrsfs,amscd,amstext,
latexsym,amsbsy,amsthm}

\usepackage{bbold}
\usepackage{mparhack}
\usepackage{marginnote}
\usepackage{xcolor}

\input{54.sty}

\nek{\toti}{%
\fontfamily{ptm}\selectfont}

\usepackage{doi}

\begin{document}


\title
{An unpublished theorem of Solovay, revisited 
}

\author 
{Ali Enayat\thanks{University of Gothenburg, 
Gothenburg, Sweden, {\tt ali.enayat@gu.se}. 
} 
\and 
Vladimir Kanovei\thanks{
IITP RAS,
 \ {\tt kanovei@googlemail.com} .
}  
\vyk{
\and 
Vassily Lyubetsky\thanks{IITP RAS,
\ {\tt lyubetsk@iitp.ru}. 
Partial support of grant RFBR 18-29-13037 acknowledged.
}} 
}

\date 
{\today}

\maketitle


\begin{abstract}
A definable pair of disjoint non-OD sets 
of reals (hence, indiscernible sets)
exists in the Sacks and $\Eo$-large 
generic extensions of 
the constructible universe $\rL$.  

{\footnotesize
\def\contentsname{\vspace*{-5ex}}
\tableofcontents
}

\end{abstract}

\punk{Introduction}
\las{int}

Let a \rit{twin partition} be any partition of a 
given set $U$ 
into two nonempty cells $A$ and $B$. 
We refer to $U$ as \rit{the universe 
of discourse}, and each of $A$ and $B$ as a twin. 
Assume that some robust notion of definability $D$ 
is chosen in advance, \eg,
$D$ might be ordinal definability \OD, 
or $D$ might be $\id11$ definability, 
or something similar.  
In this context, a twin partition 
$U = A \cup B$ can be called 
$D$-definable in one of two senses:
\bde
\item[strongly $D$-definable,]
\ie, 
each of the twins $A$ and $B$ 
is $D$-definable; 

\item[weakly $D$-definable,]
meaning that the partition $\ans{A,B}$ of $U$, 
considered as an unordered pair, is $D$-definable. 
\ede
Strong $D$-definability clearly implies weak 
$D$-definability.  
The ``twin problem'' for a given notion of definability $D$ 
is whether the converse holds. 
The twin problem obviously has a positive answer provided 
the domain of discourse $U$ contains at least one 
$D$-definable element $x$, 
then one cell of the partition consists of those  
$x^{\prime}$ that share the same cell of the partition 
as $x$, and the other cell is just the complementary set.  
This provides a trivial positive solution for the twin 
problem when $U = \omega$, or when $U$ is the class 
of ordinals,  and generally when 
$U$ admits a $D$-definable well-ordering. 
Now let's focus on the case when $U$ is the set of 
real numbers. 

The twin problem admits a positive solution 
in the case of $\id11$ definability. 
Indeed it follows from Theorem~\ref{sil} below that if 
a $\id11$ \eqr\ $\rE$ on a $\id11$ set $U$ of reals 
has precisely 
two (or even countably many) equivalence classes then 
each $\rE$-class is itself a $\id11$ set.
The problem also admits a positive solution 
in the case of $\id12$ definability because 
every non-empty $\is12$ set of reals 
contains a $\id12$ element 
(see, \eg, 4E.5 in Moschovakis~\cite{mDST}).
But slightly above of $\id12$  
there is a significant obstacle, 
as indicated by the following theorem. 

\bte
[the Sacks part originally by Solovay\snos{%
See Section~\ref{his} 
on the history of the result}]
\lam{mt}
Let\/ $a\in\dn$ be either Sacks generic 
or\/ \ela\ generic\snos
{\label{eo}%
That is, generic \poo\ the forcing by perfect sets 
$P\sq\dn$ such that the restricted relation 
${\Eo}\res P$ is not smooth, see below. 
Recall that the \eqr\ 
${\Eo}$ is defined on $\dn$ so that 
$x\Eo y$ iff the set 
$\Da(x,y)=\ens{k}{x(k)\ne y(k)}$ 
is finite.}
over\/ $\rL$.
Then it is true in\/ $\rL[a]$ that there is a\/ 
$\is12$ equivalence relation\/ $\rQ$ 
on\/ $\dn$ with exactly three equivalence 
classes, one of which is equal to\/ $\dn\cap\rL$, 
while two others are non-\OD\ sets whose union 
is equal to the\/ $\ip12$ set\/ $\dn\bez\rL$. 
\ete

Under the assumptions of this theorem, we have 
we have a \rit{weakly definable}, 
but not \rit{strongly definable}, 
partition of the $\ip12$ set $U=\dn\bez\rL$  
into two equivalence classes of $\rQ$.  
Let $A,B$ be those equivalence classes. 
As the relation $\rQ$ is lightface $\is12$, the 
unordered pair $\ans{A,B}$ is an \OD\ set, basically, 
a definable set, whose two elements 
(disjoint non-empty pointsets $A,B\sq\dn\bez\rL$) 
are non-\OD, hence, are \rit{\OD-indiscernible}. 

Models of $\zf$ or $\ZFC$ containing 
\OD\ indiscernible pairs of (non-\OD)
disjoint sets of reals are well-known.
Such is \eg\ any Sacks$\ti$Sacks extension $\rL[a,b]$
of $\rL$, where an \OD\ pair of non-\OD\ sets consists
of the $\rL$-degrees of the Sacks reals $a,b$, see
\cite{gl} and also 
\cite{ena,FGH}.
Another model with an \OD\ pair of \rit{countable}
disjoint non-\OD\ sets is defined in \cite{kl25}.
Yet those examples fail to fulfill the property 
that the union of the two sets is equal to the 
whole domain of nonconstructible reals.

Generally, \OD\ indiscernible pairs 
(not necessarily \OD\ pairs) 
of disjoint sets of reals can be extracted from  
early works on Cohen forcing. 
In particular, if $\ang{a,b}$ is a Cohen-generic, 
over $\rL$, pair of $a,b\in\dn$, then the 
$\Eo$ equivalence classes $\eko a$, $\eko b$ 
are \OD\ indiscernible in $\rL[a,b]$ 
(essentially by Feferman~\cite{fef}) and 
so are the constructibility degrees 
$[a]_\rL=\ens{x\in\dn}{\rL[x]=\rL[a]}$ 
and $[b]_\rL$~\cite{FGH}.

On the other hand, 
it is established in \cite{kl31} that, 
in some models  of $\ZFC$,
including the Sacks extension of 
the constructible universe $\rL$, it is true
that any countable
\OD\ (ordinal-definable) \rit{set of reals} consists of
\OD\ elements. 
A similar result in much more general setting is 
known from \cite[Thm 4.8]{cai} under a strong 
large cardinal hypothesis.

\punk{Outline of the proof}
\las{out}

To prove Theorem~\ref{mt}, the required \eqr\ 
will be obtained as the union of an 
increasing transfinite sequence 
$\sis{\rD_\al}{\al<\omi}$ of {\ubf countable} 
Borel \eqr s.  
The sequence is defined in $\rL$, the ground universe. 
The following is a principal definition related 
to this construction.

\bdf
\lam{dbs}
A \rit{double-bubble system}, DBS for brevity, is
a pair of \rit{countable} 
Borel \eqr s $\pae\rD\rE$
on $\dn$, such
that each $\rE$-class  
is the union of a pair of distinct $\rD$-classes.    

A DBS $\pae{\rD'}{\rE'}$ \rit{extends} $\pae\rD\rE$, 
in symbol 
$\pae\rD\rE\cle\pae{\rD'}{\rE'}$, 
if ${\rD}\sq{\rD'}$, 
${\rE}\sq{\rE'}$, and for any $x,y\in\dn,$ 
if $x\rE y$ but $x\nD y$ then we still have $x\nD' y$.
\edf

Thus the extension essentially means that the equivalence 
classes of the original \eqr s are merged in countable 
bunches, but in such a way that the two $\rD$-classes 
within the same $\rE$-class are never merged. 
We are going to define a certain $\cle$-increasing  
increasing sequence 
$\sis{\ang{\rD_\al,\rE_\al}}{\al<\omi}$ of 
double-bubble systems 
$\ang{\rD_\al,\rE_\al}$ in $\rL$, the ground universe, 
and $\rD=\bigcup_\al\rD_\al$ will be the \eqr\ 
required. 
This will take some effort.

\bpri
\lam{e0e}
The most elementary example is $\rD$= the equality, 
and $x\rE y$ iff $x(k)=y(k)$ for all $k\ge1$; $\pae\rD\rE$ 
is a DBS. 

Another example consists of the \eqr\ 
${\Eo}$ (see Footnote~\ref{eo}), 
and its subrelation $\Ee$, defined so that 
$x\Eo y$ iff the set 
$\Da(x,y)$ 
has finite even number of elements; 
$\pae\Ee\Eo$ is a DBS and obviously 
$\pae\rD\rE\cle\pae\Ee\Eo$. 
\epri

\punk{Canonization results used in the proof}
\las{ing}

Here we present some well-known results of 
modern descriptive set theory involved in the 
proof of Theorem~\ref{mt}.
We begin with the Silver Dichotomy theorem 
and a canonization corollary. 
See \eg\ \cite[Theorem 2.2]{sami19} or 
\cite[Section 10.1]{kB} for a 
proof of the ``moreover'' lightface version 
of Theorem~\ref{sil}. 

\bte
[Silver's Dichotomy \cite{sil}]
\lam{sil}
Suppose that\/ $\rE$ is a\/ $\fp11$ \eqr\ on 
a Borel set\/ $X\sq\dn.$ 
Then either\/ $\rE$ has at most countably many 
equivalence classes, or there exists  
a perfect partial\/ $\rE$-transversal\snos
{A \rit{partial transversal} is a set of pairwise 
inequivalent elements. 
A \rit{full} transversal requires that in addition 
it has a non-empty intersection with any 
equivalence class in a given domain.}. 

If moreover\/ $X$ is lightface\/ $\id11$ and\/ 
$\rE$ is lightface\/ $\ip11$ then all equivalence 
classes are lightface\/ $\id11$ 
in the ``either'' case.\qed
\ete


\bcor
\lam{silC}
Suppose that\/ $\rE$ is a\/ $\fp11$ \eqr\ on 
a Borel set\/ $X\sq\dn.$ 
Then there is a perfect set\/ $Y\sq X$ such that\/ 
$\rE$ coincides on\/ $Y$ with$:$ 
\bit
\item[$-$] 
either\/ {\rm(I)} the total equivalence $\tot$ 
making all reals equivalent$;$

\item[$-$] 
or\/ {\rm(II)} 
the equality, so that\/ $Y$ is a partial\/ 
$\rE$-transversal. 
\eit
If in addition\/ $\rE$ is countable\snos
{An \eqr\ is \rit{countable} iff all its equivalence 
classes are at most countable.}
then\/ {\rm(I)} is impossible. 
\ecor
\bpf
In the ``or'' case of Theorem~\ref{sil} we have (II). 
In the ``either'' case pick an uncountable 
equivalence class $C$ and let $Y\sq C$ be any perfect 
set.
\epf 

\bcor
\lam{silF}
If\/ $X\sq\dn$ is a perfect set, and\/ 
$f:X\to\dn$ a \kb\ map,  
then there is a perfect set\/ 
$Y\sq X$ such that\/ 
$f\res Y$ is a bijection or a constant.
\ecor
\bpf
This is a well-known fact, of course, yet it    
immediately follows from Corollary~\ref{silC}. 
Indeed define a Borel \eqr\ $\rE$ on $X$ 
such that $x\rE y$ iff $f(x)=f(y)$. 
Apply Corollary~\ref{silC}. 
\epf

\gol{Now we recall some definitions and 
results related to} \ela\ sets.
A Borel set $X\sq\dn$ is called \rit\ela\ if 
${\Eo}\res X$ is still a non-smooth\snos
{\label{smu}%
Recall that an \eqr\ $\rE$ on a Borel set $X$ 
is \rit{smooth} if there is a Borel map $f:X\to\dn$ 
such that we have $x\rE y$ iff $f(x)=f(y)$ for 
all $x,y\in X$. 
The \eqr\ $\Eo$ is non-smooth on $\dn,$ 
meaning that such a Borel $f$ does not exist. 
See Example 6.5 in \cite{kemi}.} 
\eqr. 
For instance $\dn$ itself is \ela, while 
any Borel partial $\Eo$-transversal is not. 
If $\bu=\sis{u_n^i}{n<\om, i=0,1}$ is an 
array of strings $u_n^i\in\bse,$ 
satisfying $\lh{u_n^0}=\lh{u_n^1}\ge1$ and 
$u_n^0\ne u_n^1$ for all $n$, 
then we call $\bu$ \gol{\rit{a \ema}}, let 
$$
x_\bu^a=u_0^{a(0)}\we u_1^{a(1)}\we u_2^{a(2)}\we
\dots\we u_n^{a(n)}\we\dots \in \dn\,.
$$
for any $a\in\dn$, and 
define a \rit\cela\ set 
$\dxu\bu=\ens{x_\bu^a}{a\in\dn}$. 
Each \cela\ set $\dxu\bu$ is perfect, 
and \ela\ via the map 
$a\mto x_\bu^a$. 
On the other hand, 
it is known (see \eg\ \cite[Section 7.1]{ksz}) 
that  each (Borel) \ela\ set $X\sq\dn$ contains a 
\cela\ subset $Y\sq X$. 

\gol{
If further $\bv=\sis{v_n^i}{n<\om, i=0,1}$ is another 
\ema, then we define a   
homeomorphism and $\Eo$-isomorphism 
$h_{\bu\bv}:\dxu\bu\onto \dxu\bv$ 
such that $h_{\bu\bv}(x_\bu^a)=x_\bv^a$ 
for all $a\in\dn$.
Maps of the form $h_{\bu\bv}$ will be called 
\rit{\cela\ maps}.
}

\bte
[Theorem 7.1 in \cite{ksz}, or else \cite{millBcanon}]
\lam{ksz1}
Suppose that\/ $\rE$ is a Borel \eqr\ on\/ $\dn,$ 
and\/ $X\sq\dn$ is a\/ \ela\ set. 
Then there is a \cela\ set\/ $Y\sq X$ such that\/ 
$\rE$ coincides on\/ $Y$ with$:$
\bit
\item[$-$] 
either\/ {\rm(I)} the total \eqr\/ $\tot$$;$

\item[$-$] 
or\/ {\rm(II)} the relation\/ $\Eo\;;$

\item[$-$] 
or\/ {\rm(III)} 
 the equality. 
\eit
In addition, if\/ $\rE$ is a countable \eqr\ then\/ 
{\rm(I)} is impossible, while if\/ ${\Eo}\sq{\rE}$ 
then\/ {\rm(III)} is impossible.\qed 
\ete

\bcor
\lam{kszC}
If\/ $X\sq\dn$ is a Borel\/ \ela\ set, and\/ 
$Z\sq X$ a Borel set,  
then there is a \cela\ set 
$Y\sq X$ such that\/ 
$Y\sq Z$ or\/ $Y\cap Z=\pu$.
\ecor
\bpf
Define a Borel \eqr\ $\rE$ on $X$ such that 
$x\rF y$ iff $x,y\in Z$ or $x,y\in X\bez Z$. 
Apply Theorem~\ref{ksz1}. 
As $\rE$ has just two equivalence classes, 
only (I) is possible. 
\epf

\bcor
\lam{kszF}
If\/ $X\sq\dn$ is a Borel\/ \ela\ set, and\/ 
$f:X\to\dn$ a \kb\ map,  
then there exists a \cela\ set\/ 
$Y\sq X$ such that\/ 
$f\res Y$ is a bijection or a constant.
\ecor
\bpf
Define a Borel \eqr\ $\rE$ on $X$ such that 
$x\rE y$ iff $f(x)=f(y)$. 
Apply Theorem~\ref{ksz1}. 
We have to prove that (II) is impossible. 
Suppose to the contrary that ${\rE}={\Eo}$ on 
a \cela\ set $Y\sq X$. 
In other words, we have $f(x)=f(y)$ iff $x\Eo y$ 
for all $x,y\in Y$. 
Thus $f$ is a \kb\ reduction of ${\Eo}\res Y$ to 
the equality, which contradicts to the assumption 
that $Y$ is \ela.
\epf

As a forcing notion, the set 
$\peo$ of all \cela\ (perfect) sets adjoins    
reals of minimal degree, preserves $\ali$, and 
has some other remarkable properties 
resembling the Sacks forcing, 
see \eg\ \cite[Section 7.1]{ksz} 
and references thereof.

\punk{\Cor ing maps, Sacks case}
\las{corS}

 \bdf
\lam{kap}
Given a set $X\sq\dn$ and a map $f:X\to\dn$, a
DBS $\pae\rD\rE$:\vom

--
\rit{\cor s} $f$
if $f(x)\in\eke x$ for all $x\in X$;\vom

-- 
\rit{positively \cor s} $f$
if $f(x)\in\ekd x$ for all $x\in X$;\vom

--  
\rit{negatively \cor s} $f$
if $f(x)\in\eke x\bez\ekd x$ for all $x\in X$. 
\edf

\ble
\lam{s2}
Assume that\/ $\pae\rD\rE$ is a DBS, 
$X\sq\dn$ is a perfect set, 
and\/ $f:X\to\dn$ is \kb\  and 1-1.
There exist a\/ perfect set\/ 
$Y\sq X$ and a DBS\/ $\pae{\rD'}{\rE'}$ 
which extends\/ $\pae\rD\rE$ and   
\cor s\/ $f\res Y$.
\ele
\bpf
The sets $X'=\ens{x\in X}{x\rE f(x)}$ and
$X''=\ens{x\in X}{x\nE f(x)}$ are Borel, hence
there is a perfect set $X_0$ with either $X_0\sq X'$
or $X_0\sq X''$.
But if $X_0\sq X'$ then $\pae\rD\rE$ already
\cor s\/ $f\res X_0$, and we are done.
Thus we assume that $X_0\sq X''$, that is, $x\nE f(x)$
for all $x\in X_0$.

By Corollary~\ref{silC}, there is a perfect 
set $X_1\sq X_0$ such that $\rE,\rD$   
coincide with the equality on $X_1$. 
Define an \eqr\ $\wE$ on $X_1$ such that 
$x\wE y$ iff $f(x)\rE f(y)$, 
and define $\wD$ similarly. 
Consider the $\sq$-minimal \eqr\ $\rF$ defined 
on $\dn$ such that ${\rE}\sq{\rF}$ and if 
$x,y\in \dn$ and $f(x)\rE y$ then 
$x\rF y$. 
Thus $\wE,\wD,\rF$ are countable Borel 
\eqr s on $X_1$. 
(The borelness of $\rF$ 
holds since all intended quantifiers 
in the definition of $\rF$ are over countable 
domains.) 
By Corollary~\ref{silC}, there is a
perfect set $Y\sq X_1$ such that $\wE,\wD,\rF$
coincide with the equality on $Y$, along with $\rE,\rD$.
It follows, by the choice of $X_0$, that if
$x,y\in Y$ (whether equal or not) then $x\nE f(y)$.


We define the \eqr s $\rE',\rD'$ as follows. 

If $x\in\dn$ 
and the $\rE$-class $\eke x$ 
does {\ubf not} intersect 
\rit{the critical domain} 
$\Da=Y\cup\ens{f(x)}{x\in Y},$ 
then put $\ekep x=\eke x$ and 
$\ekdp x=\ekd x$, so such a $\rE$-class 
and its $\rD$-subclasses are not changed. 
But within $\Da$ some classes will be merged. 
Namely if $x\in Y$ then we have to merge 
$\eke x$ with $\eke {f(x)}$, hence put
$$
\ekep x
=
\eke{x}\cup
\eke{f(x)} 
\qand
\ekdp x
=
\ekd{x}\cup\ekd{f(x)}\,,
$$
and define the other $\rD'$-class within 
$\ekep x$ as $\ekep x\bez\ekdp x$.

A routine verification shows that in either  
case the relations $\rE',\rD'$ are Borel, and 
the pair $\pae{\rD'}{\rE'}$ is a DBS 
which extends $\pae{\rD}{\rE}$ and 
positively \cor s $f\res Y$ 
(because we have $f(x)\in\ekdp x$ 
for all $x\in Y$ simply by construction).
\epf

\ble
\lam{s3}
Let\/ $\pae\rD\rE$ be a DBS, and\/
$R,X\sq\dn$ be perfect sets.
There exist$:$ 
a\/ perfect set\/ $Y\sq X,$ 
\kb\  1-1 maps\/ $f,g:Y\to R,$ 
and a DBS\/ $\pae{\rD'}{\rE'}$ which extends\/ 
$\pae\rD\rE$, 
positively \cor s~$f\res Y,$  
and negatively \cor s~$g\res Y.$ 
\ele
\bpf
By Corollary~\ref{silC}, there exist perfect 
partial $\rE$-transversals $X_0\sq X$ 
and $R_0\sq R$. 
Let $R_0=R_1\cup R_2$ be a partition into two 
disjoint perfect sets.
Then $\eke{R_1}$ and $\eke{R_2}$ are disjoint,
hence there is a perfect set $Y\sq X_0$ such
that $\eke{Y}$ does not intersect either 
$\eke{R_1}$ or $\eke{R_2}$.
Let say $\eke{Y}\cap\eke{R_1}=\pu$.

Let $R_1=R'\cup R''$ be a partition into two 
disjoint perfect sets.
It follows by construction that (*)
the Borel sets\/ $Y,R',R''$ 
are pairwise disjoint and the union\/ 
$\Da=Y\cup R'\cup R''$ is a 
partial\/ $\rE$-transversal.
Let $f:Y\to R'$ and $g:Y\to R''$ be 
arbitrary \kb\  1-1 maps.  

We define the \eqr s $\rE',\rD'$ as follows. 

If $x\in\dn$ 
and the $\rE$-class $\eke x$ 
does {\ubf not} intersect 
\rit{the critical domain} 
$\Da=Y\cup Z'\cup Z'',$ 
then put $\ekep x=\eke x$ and 
$\ekdp x=\ekd x$, so such a $\rE$-class 
and its $\rD$-subclasses are not changed. 
But within $\Da$ some classes will be merged. 
Namely if $x\in Y$ then we have to merge 
$\eke x$ with $\eke {f(x)}$ and $\eke {g(x)}$, 
hence we put 
$\ekep x=\eke{x}\cup\eke{f(x)}\cup\eke{g(x)}$. 
We further define  
$$
\ekdp x=\ekd{x}\cup\ekd{f(x)}\cup
{(\eke{g(x)}\bez\ekd{g(x)})}\,,
$$
and let 
${(\eke{x}\bez\ekd{x})}\cup
{(\eke{f(x)}\bez\ekd{f(x)})}\cup\ekd{g(x)}$
be the other $\rD'$-class within 
$\ekep x$.
A routine verification using (*)  
shows that the relations $\rE',\rD'$ are Borel, 
and the pair $\pae{\rD'}{\rE'}$ is a DBS 
that extends $\pae{\rD}{\rE}$, 
positively \cor s $f\res Y$, and 
negatively \cor s $g\res Y$.
\epf

\punk{\Cor ing maps, \ela\ case}
\las{corE}

Here we prove two \cor ing lemmas similar 
to \ref{s2} and \ref{s3}, yet with somewhat 
more complex proofs.

\ble
\lam{e2}
Assume that\/ $\pae\rD\rE$ is a DBS, 
${\Eo}\sq{\rE}$, 
$X\sq\dn$ is a \cela\ set, 
and\/ $f:X\to\dn$ is \kb\  and 1-1.
There exist a\/ \cela\ set 
$Y\sq X$ and a DBS\/ $\pae{\rD'}{\rE'}$ 
which extends\/ $\pae\rD\rE$ and 
\cor s\/ $f\res Y$.
\ele
\bpf
First of all, arguing as in the proof of
Lemma~\ref{s2} (but using Corollary~\ref{kszC}),
we get a \cela\ set $X_0\sq X$ with
$x\nE f(x)$ for all $x\in X_0$.
By Theorem~\ref{ksz1}, there is a \cela\ perfect 
set $X_1\sq X_0$ such that the relations $\rE,\rD$   
coincide with $\Eo$ on $X_1$. 
Define an \eqr\ $\wE$ on $X_1$ such that 
$x\wE y$ iff $f(x)\rE f(y)$, 
and define $\wD$ similarly. 
Consider the $\sq$-minimal \eqr\ $\rF$ defined 
on $\dn$ such that ${\rE}\sq{\rF}$ and if 
$x,y\in\dn$ and $f(x)\rE y$ then 
$x\rF y$. 
Thus $\wE,\wD,\rF$ are countable Borel 
\eqr s on $X_1$. 
(The borelness of $\rF$ 
holds since all intended quantifiers 
in the definition of $\rF$ are over countable 
domains.) 
By Theorem~\ref{ksz1}, there is a \cela\ perfect 
set $Y\sq X_1$ such that each of these three 
\eqr s is either of type (I) or of type (II) 
on $Y$. 
However, as each $\rE$-class contains two 
$\rD$-classes, $\wE$ has to coincide with $\wD$ on 
$Y$. 
Finally, as ${\rE}\sq{\rF}$, we have
${\rF}={\Eo}$ on $Y.$
It follows by the choice of $X_0$ that if $x,y\in Y$
(whether equal or not) then $x\nE f(y)$.

To conclude, ${\rE}={\rD}={\rF}={\Eo}$ on $Y$, 
and also either ${\wE}={\wD}$ is the equality 
on $Y$, or ${\wE}={\wD}={\Eo}$ on $Y$.
This leads to the following two cases. 

In each case, we are going to define the 
\eqr s $\rE',\rD'$ required. 
If $x\in\dn$ 
and the $\rE$-class $\eke x$ 
does {\ubf not} intersect 
\rit{the critical domain} 
$\Da=Y\cup\ens{f(x)}{x\in Y},$ 
then put $\ekep x=\eke x$ and 
$\ekdp x=\ekd x$, so such a $\rE$-class 
and its $\rD$-subclasses are not changed. 
But within $\Da$ some classes will be merged. 
In particular, we are going to merge 
$\eke x$ with $\eke {f(x)}$ 
for any $x\in Y$.\vom 

{\ubf Case 1:} ${\wE}={\wD}$ is the equality on 
$Y$ while ${\rD}={\rE}={\rF}={\Eo}$ on $Y$, 
thus if $x, y\in Y$ then first, 
$x\ne y$ implies $f(x)\nE f(y)$ and 
$f(x)\nD f(y)$, 
and second, 
$ 
\eke x\cap Y=\ekd x\cap Y=\eko x\cap Y\,. 
$
If $x\in Y$ then put
$$
\ekep x
=
\eke{x}\cup
{\bigcup_{y\in Y\cap\eko x}\eke{f(y)}}
\qand
\ekdp x
=
\ekd{x}\cup
{\bigcup_{y\in Y\cap\eko x}
\ekd{f(y)}}\,,
$$
and define the other $\rD'$-class within 
$\ekep x$ as $\ekep x\bez\ekdp x$.\vom

{\ubf Case 2:} 
${\rE}={\rD}={\wE}={\wD}={\rF}={\Eo}$
on $Y$, that is, if $x, y\in Y$ then  
$$ 
{x\Eo y}\leqv {x\rE y}\leqv {x\rD y} 
\leqv {f(x)\rE f(y)}\leqv {f(x)\rD f(y)}\,. 
$$ 
Assume that $x\in Y.$ 
Put    
$\ekep x=\eke{x}\cup\eke{f(x)}=
\eke{y}\cup\eke{f(y)}$ 
for any other $y\in Y\cap\eko x$, and 
$\ekdp x=\ekd{x}\cup\ekd{f(x)}=
\ekd{y}\cup\ekd{f(y)}$ 
for any other $y\in Y\cap\eko x$.
Define the other $\rD'$-class within 
$\ekep x$ as $\ekep x\bez\ekdp x$. 

A routine verification shows that in either  
case the relations $\rE',\rD'$ are Borel, and 
the pair $\pae{\rD'}{\rE'}$ is a DBS 
which extends $\pae{\rD}{\rE}$ and 
\cor s $f\res Y$ 
(because we have $f(x)\in\ekep x$ 
for all $x\in Y$ simply by construction).
\epf

\ble
\lam{e3}
Let\/ $\pae\rD\rE$ be a DBS with\/ ${\Eo}\sq{\rE}$,
and\/ $R,X\sq\dn$ be\/ \cela\ sets.
There exist$:$ 
a\/ \cela\ set\/ $Y\sq X,$ 
\gol{\cela\  maps\/} $f,g:Y\to R,$ 
and a DBS\/ $\pae{\rD'}{\rE'}$ that extends\/ 
$\pae\rD\rE$, positively \cor s~$f,$  
and negatively \cor s~$g.$ 
\ele
\bpf
\gol{%
By Theorem~\ref{ksz1}, we \noo\ assume that $\rE$ 
coincides with $\Eo$ on $R$. 
By definition, $R=\dxu\br$ for a \ema\ 
$\br=\sis{r^i_n}{n<\om,i=0,1}$. 
Now let 
$\bp=\sis{p^i_n}{n<\om,i=0,1}$, 
$\bq=\sis{q^i_n}{n<\om,i=0,1}$, 
where 
$p^i_n=r^0_{2n}\we r^i_{2n+1}$,  
$q^i_n=r^1_{2n}\we r^i_{2n+1}$. 
Thus $\bp,\bq$ are \emas, and the sets 
$\dxu\bp$, $\dxu{\bq}$ satisfy 
$\dxu\bp\cup\dxu{\bq}\sq\dxu\br=R$ and 
$\eko{\dxu\bp}\cap\eko{\dxu{\bq}}=\pu$, 
hence, 
$\eke{\dxu\bp}\cap\eke{\dxu{\bq}}=\pu$ by the 
assumption above. 
It follows by Corollary~\ref{kszC} that there is a 
\cela\ set $X_0\sq X$ satisfying
$\eke{X_0}\cap\eke{\dxu\bp}=\pu$ or 
$\eke{X_0}\cap\eke{\dxu{\bq}}=\pu$. 
Let say $\eke{X_0}\cap\eke{\dxu\bp}=\pu$. 
As just above, there exist \emas\ 
$\bp',\bp''$ such that the \cela\ sets 
$R'=\dxu{\bp'}$, $R''=\dxu{\bp''}$ satisfy 
$R'\cup R''\sq\dxu\bp$ and 
$\eke{R'}\cap\eke{R''}=\pu$. 

To conclude, we have \cela\ sets $X_0\sq X$ and 
$R',R''\sq R$ satisfying 
$\eke{R'}\cap\eke{R''}=
\eke{X_0}\cap\eke{R'}=
\eke{X_0}\cap\eke{R''}=\pu$.%
}
Theorem~\ref{ksz1} yields a 
\cela\ set $Y=\dxu\bu\sq X_0$ such that 
${\rE}={\rD}={\Eo}$ on $Y$. 
\gol{%
Consider the \cela\ maps $f=h_{\bu\bp'}:Y\to R'$ and 
$g=h_{\bu\bp''}:Y\to R''$.
}

We define the \eqr s $\rE',\rD'$ as follows. 

If $x\in\dn$ 
and the $\rE$-class $\eke x$ 
does {\ubf not} intersect 
\rit{the critical domain} 
$\Da=Y\cup (\imx f{Y})\cup (\imx g{Y}),$ 
then put $\ekep x=\eke x$ and 
$\ekdp x=\ekd x$, so such a $\rE$-class 
and its $\rD$-subclasses are not changed. 
But within $\Da$, if $x\in Y$ then we have to merge 
$\eke x$ with $\eke {f(x)}$ and $\eke {g(x)}$, 
hence we put
$$
\ekep x
=
\eke{x}\cup
{\bigcup_{y\in Y\cap\eko x}\eke{f(y)}}
\qand 
\ekdp x
=
\ekd{x}\cup
{\bigcup_{y\in Y\cap\eko x}
\ekd{f(y)}}\,,
$$
and define the other $\rD'$-class within 
$\ekep x$ as $\ekep x\bez\ekdp x$.
\vyk{
$\ekep x=\eke{x}\cup\eke{f(x)}\cup\eke{g(x)}$. 
We further define  
$$
\ekdp x=\ekd{x}\cup\ekd{f(x)}\cup
{(\eke{g(x)}\bez\ekd{g(x)})}\,,
$$
and let 
${(\eke{x}\bez\ekd{x})}\cup
{(\eke{f(x)}\bez\ekd{f(x)})}\cup\ekd{g(x)}$
be the other $\rD'$-class within 
$\ekep x$.
}%
A routine verification
shows that the relations $\rE',\rD'$ are Borel, 
and the pair $\pae{\rD'}{\rE'}$ is a DBS 
that extends $\pae{\rD}{\rE}$, 
positively \cor s $f\res Y$, and 
negatively \cor s $g\res Y$.
\epf

\punk{Increasing system of \eqr s}
\las{ise}

\bpro
[{\ubf in $\rL$}]
\lam{iseS}
There is an\/ $\cle$-increasing sequence 
of DBSs\/ $\pae{\rD_\al}{\rE_\al}$, $\al<\omi$, 
beginning with\/ $\Eo$ of Footnote~\ref{eo} 
and\/ ${\rD_0}={\Ee}$ and such that 
\ben
\renu
\itlb{iseS1}%
if\/ $X\sq\dn$ is perfect and\/ $f:X\to\dn$ \kb\  
and 1-1, then there exist$:$ 
a perfect\/ $X'\sq X$ and 
an ordinal\/ $\al<\omi$ such that\/  
$\pae{\rD_\al}{\rE_\al}$ \cor s\/ $f\res X'$$;$ 

\itlb{iseS2}%
if\/ $X,R\sq\dn$ are perfect sets, 
then there exist$:$ 
a perfect set\/ $Y\sq X$,  
an ordinal $\al<\omi$,  
and \kb\  1-1 maps\/ $f,g:Y\to R$,  
such that\/  
$\pae{\rD_\al}{\rE_\al}$ \cor s\/ $f$ positively 
and \cor s\/ $g$ negatively$;$

\itlb{iseS3}%
the sequence of pairs\/ $\pae{\rD_\al}{\rE_\al}$ 
is\/ $\id12$, in the sense that there exists 
a\/ $\id12$ sequence of codes for Borel sets\/ 
$\rD_\al$ and\/ $\rE_\al$. 
\een
\epro
\bpf
An obvious inductive construction using lemmas 
\ref{s2}, \ref{s3}, that takes 
a G\"odel-least code of all possible pairs fitting 
the given inductive step, with the obvious 
union at limit steps . 
\epf

\bpro
[{{\ubf in $\rL$}}]
\lam{iseE}
There is an\/ $\cle$-increasing sequence 
of DBSs\/ $\pae{\rD_\al}{\rE_\al}$, $\al<\omi$, 
beginning with\/ $\Eo$ of Footnote~\ref{eo} 
and\/ ${\rD_0}={\Ee}$ and such that 
\ben
\renu
\itlb{iseE1}%
if\/ $X\sq\dn$ is a Borel\/ \ela\ set 
and\/ $f:X\to\dn$ 
\kb\  
and 1-1, then there exist$:$ 
a\/ \cela\ set\/ $Y\sq X$ 
and\/ $\al<\omi$ such that\/  
$\pae{\rD_\al}{\rE_\al}$ \cor s\/ $f\res Y$$;$ 

\itlb{iseE2}%
\gol{if\/ $X,R\sq\dn$ are\/ \ela\ sets, 
then there exist$:$  
a\/ \cela\/ set\/ $Y\sq X$,  
an ordinal $\al<\omi$, 
and \cela\ maps\/ $f,g:Y\to R$,  
such that\/  
$\pae{\rD_\al}{\rE_\al}$ \cor s\/ $f$ positively 
and\/ $g$ negatively$;$}

\itlb{iseE3}%
the sequence of pairs\/ $\pae{\rD_\al}{\rE_\al}$ 
is\/ $\id12$, in the sense that there exists 
a\/ $\id12$ sequence of codes for Borel sets\/ 
$\rD_\al$ and\/ $\rE_\al$. 
\een
\epro
\bpf
Similar. 
\epf

\punk{Proof of the main theorem}
\las{mat}

\bpf[Theorem~\ref{mt}, Sacks case]
Fix, in $\rL$, an $\cle$-increasing sequence 
of DBSs $\pae{\rD_\al}{\rE_\al}$, $\al<\omi$, 
satisfying conditions 
\ref{iseS1}, \ref{iseS2}, \ref{iseS3} 
of Proposition~\ref{iseS}. 

{\ubf Arguing in a Sacks-generic extension 
$\rL[a_0]$,} 
we define a relation $\rD=\bigcup_{\al<\omi}\rD_\al$ 
on $\dn$; thus
$x\rD y$ iff $x\rD_\al y$ for some $\al<\omi$. 
(We identify Borel sets $\rD_\al$ and $\rE_\al$, 
formally defined in $\rL$, with their extensions,  
Borel sets in $\rL[a_0]$ with the same codes.) 
Define a relation $\rE=\bigcup_{\al<\omi}\rE_\al$ 
 on $\dn$ similarly. 
Define the subdomain\/ $U=\dn\bez\rL$ 
of all new reals. 
Then $a_0\in U$ and  
all reals in $U$ have the same $\rL$-degree by the 
minimality of Sacks reals, see \eg\  
\cite[Theorem 15.34]{jechmill}.

\ble
\lam{les}
It is true in\/ $\rL[a_0]$ that 
\ben
\renu
\itlb{les1}%
$\rE$ and\/ $\rD$ are \eqr s and\/ $\rD$ is 
a subrelation of\/ $\rE\;;$

\itlb{les0}%
$\rD$ is lighface\/ $\is12\;;$

\itlb{les2}%
all reals\/ $x,y\in U$ are\/ $\rE$-equivalent$;$ 

\itlb{les3}%
there are exactly two\/ $\rD$-classes 
intersecting\/ $U$ --- call them\/ $M\yi N\;;$  

\itlb{les4}%
the sets\/
$M,N$ are not\/ \OD\snos
{Note that $M,N$ are indiscernible 
in a stronger sense: if $R(M,N)$ holds for some 
\OD\ relation $R$, then $R(N,M)$ holds. 
Indeed, otherwise $M$ can be distinguished from 
$N$ by the property: 
``$R(\cdot,A)$ holds but $R(A,\cdot)$ fails, where $A$ 
is the other element of the pair $\ans{M,N}$''.}, 
hence\/ $M\cup N=U$. 
\een
\ele
\bpf
\ref{les1}
To see that $\rE$ is an \eqr, let $a,b,c\in W$ 
and suppose that $a\rE b$ and $a\rE c$. 
Then by definition we have $a\rE_\al b$ 
and $a\rE_\al c$ for some $\al<\omi$. 
However being an \eqr\ is absolute by 
Shoenfield's absoluteness theorem 
\cite[Theorem 25.20]{jechmill}. 
Therefore $b\rD_\al c$ holds, as required.

\ref{les0} 
holds by 
Theorem~\ref{iseS}\ref{iseS3}. 

\ref{les2} 
Let $b\in U$; 
\gol{prove that $a_0\rE b$. 
It is a known property of the Sacks forcing 
that there is a \kb\ 1-1 map $f:\dn\to\dn$ 
with a code in $\rL$, such that 
$b=f(a_0)$.\snos
{\label{bico}%
\gol{Indeed, 
by the property of \kb\  reading of names, 
we have $b=f(a)$, where 
$f:\dn\to\dn$ is a \kb\  map with a code 
in $\rL$. 
But any \kb\  $g:X\to\dn,$ 
defined on a perfect set $X\sq\dn,$ is 1-1 
or a constant on a smaller perfect set, by 
Corollary~\ref{silF}.  
Thus   
there is a perfect set $Y\sq\dn,$
coded in $\rL$,
such that $a_0\in Y$ and $f\res Y$ is 1-1 
or a constant. 
However if $f$ is a constant, say $f(x)=z_0\in\dn$ 
for all $x\in Y$, then $f(a_0)=b=z_0\in\rL$, 
which contradicts to $b\nin\rL$.}
} 
It follows then from 
Theorem~\ref{iseS}\ref{iseS1} that there exists 
a perfect set $X\sq\dn,$ coded in $\rL$ and 
such that $a_0\in X$ and $\rE_\al$ \cor s $f\res X$ 
for some $\al$. 
In particular, $\ang{a_0,b}\in{\rE_\al}$, hence 
we have $a_0\rE b$ as required.} 

\ref{les3} 
Let $a,b,c\in U$; 
prove that two of these reals are 
$\rD$-equivalent. 
Note that $a\rE b\rE c$ by \ref{les2}, 
and hence there is an ordinal $\al<\omi$ 
such that $a\rE_\al b\rE_\al c$. 
However containing exactly two $\rD_\al$-classes in 
each $\rE_\al$-class is absolute. 
It follows that at least one pair among $a,b,c$ 
is $\rD_\al$-equivalent, as required. 

\ref{les4} 
Suppose to the contrary that $M$ is \OD. 
Then $M$ is Sacks-forced over $\rL$, meaning that 
there is a perfect set $R\sq\dn,$ coded in $\rL$ 
and such that $R\cap U\sq M$ in $\rL[a_0]$. 
By Proposition~\ref{iseS}\ref{iseS2}, there exist: 
a perfect set $Y\sq \dn$ 
\gol{coded in $\rL$ and} containing $a_0$,  
an ordinal $\al<\omi$,  
and \kb\  1-1 maps $f,g:Y\to R$, 
\gol{also coded in $\rL$ and}
such that  
${\rE_\al}$ \cor s $f\res Y$ positively 
and $g\res Y$ negatively.
In other words the reals $b=f(a_0)$ 
and $c=g(a_0)$ in $U\cap R$ satisfy 
$a_0\rD_\al b$, $a_0\rE_\al c$, 
but $\neg\: (a_0\rD_\al c)$. 
It easily follows that $b\nD c$, which 
contradicts the fact that $b,c$ belong to 
one and the same $\rD$-class.
\epf

To conclude, it is true in the Sacks extension 
$\rL[a_0]$ that $\rD$ is a $\is12$ \eqr\ on $\dn,$ 
and the nonconstructible domain $U=\dn\bez\rL$ 
(a $\ip12$ set) is equal to 
the union of two (non-empty) $\rD$-equivalence 
classes, which are non-\OD\ sets.
Now, to prove Theorem~\ref{mt} (Sacks case),
it suffices to define the required \eqr\ $\rQ$ 
on $\dn$ in $\rL[a_0]$ as follows:
$x\rQ y$ iff $x\rD y$ 
or just $x,y$ both belong to $\rL$. 
\epF{Theorem~\ref{mt}, Sacks case}

\bpf[Theorem~\ref{mt}, \ela\ case]
Rather similar to the proof of 
in the Sacks case above. 
{\ubf Arguing in a \ela-generic extension 
$\rL[a_0]$,} 
we define relations $\rD=\bigcup_{\al<\omi}\rD_\al$, 
$\rE=\bigcup_{\al<\omi}\rE_\al$  on $\dn,$ and
the subdomain\/ $U=\dn\bez\rL$; $a_0\in U$.

\gol{
\ble
\lam{lee}
It is true in\/ $\rL[a_0]$ that 
\ben
\renu
\itlb{lee1}%
$\rE$ and\/ $\rD$ are \eqr s and\/ $\rD$ is 
a subrelation of\/ $\rE\;;$

\itlb{lee0}%
$\rD$ is lighface\/ $\is12\;;$

\itlb{lee2}%
all reals\/ $x,y\in U$ are\/ $\rE$-equivalent$;$ 

\itlb{lee3}%
there are exactly two\/ $\rD$-classes 
intersecting\/ $U$ --- call them\/ $M\yi N\;;$  

\itlb{lee4}%
the sets\/
$M,N$ are not\/ \OD, 
hence\/ $M\cup N=U$. 
\een
\ele
\bpf
The proof of claims 
\ref{lee1}, \ref{lee0}, \ref{lee2}, \ref{lee3} 
goes on similarly to Lemma~\ref{les},  
with some obvious changes 
{\sl mutatis mutandis\/}, in particular, the 
reference to Corollary~\ref{silF} has to be 
replaced by Corollary~\ref{kszF} in Footnote~\ref{bico}, 
the Proposition~\ref{iseS} 
by Proposition~\ref{iseE}, and so on. 
But the last claim needs special attention because 
not all new reals in $\rL[a_0]$ are \ela-generic 
unlike the Sacks case.

\ref{lee4}
First of all let's prove that each of the classes $M,N$ 
of \ref{lee3} contains a real $b\in\dn$ \ela-generic 
over $\rL$. 
Indeed in view of \ref{lee3} it suffices to prove that 
(*) 
\rit{there are\/ \ela-generic, but not\/ $\rD$-equivalent, 
reals\/ $b,c\in\rL[a_0]\cap\dn.$} 
Emulating the proof of Theorem~\ref{les}\ref{les4}, 
but using \ref{iseE}\ref{iseE2}
instead of~\ref{iseS}\ref{iseS2}, we find a  
\cela\ set $Y\sq \dn,$ coded in $\rL$ and containing $a_0$,  
an ordinal $\al<\omi$,  
and \cela\ maps $f,g:Y\to \dn$, also coded in $\rL$ and  
such that  
${\rE_\al}$ \cor s $f\res Y$ positively 
and $g\res Y$ negatively.
We conclude that the reals $b=f(a_0)$ 
and $c=g(a_0)$ in $U$ satisfy 
$a_0\rD_\al b$, $a_0\rE_\al c$, 
but $\neg\: (a_0\rD_\al c)$, so that $b\nD c$. 
And finally, it is clear that $b,c$ are \ela-generic 
along with $a_0$. 
(Basically any image of a \ela-generic real $a\in\dn$ 
via a \cela\ map $h$, coded in $\rL$, 
with $a\in\dom h$, is \ela-generic 
by an easy argument.)

Now suppose to the contrary that $M$ is \OD. 
Let $\mu(\cdot)$ be an $\in$-formula, with ordinals 
as parameters, such that $M=\ens{x}{\mu(x)}$ in $\rL[a_0]$. 
By (*), there is a real $b_0\in M$ (in $\rL[a_0]$), 
\ela-generic over $\rL$. 
Then it is true in $\rL[a_0]=\rL[b_0]$ that $\mu(b_0)$ 
and any real $x$ satisfying $\mu(x)$ also satisfies $x\rD b_0$. 
This is \ela-forced over $\rL$, meaning that 
there is a \cela\ set $R\sq\dn,$ coded in $\rL$ 
and such that (1) $b_0\in R$, 
(2) every real $b\in R\cap\rL[x]$, \ela-generic over $\rL$, 
satisfies $\mu(x)$ in $\rL[b]=\rL[b_0]=\rL[a_0]$, 
and hence satisfies $b\rD b_0$. 

However, emulating the proof of Theorem~\ref{les}\ref{les4} 
as above, we find a  
\cela\ set $Y\sq \dn,$ coded in $\rL$ and containing $b_0$,  
an ordinal $\al<\omi$,  
and \cela\ maps $f,g:Y\to R$, also coded in $\rL$ and  
such that  
${\rE_\al}$ \cor s $f\res Y$ positively 
and $g\res Y$ negatively.
Then the reals $b=f(b_0)$ and $c=g(b_0)$ are 
\ela-generic over $\rL$ and satisfy 
$b_0\rD_\al b$ and $b_0\rE_\al c$ but 
$\neg\: (b_0\rD_\al c)$, hence  $b\nD c$, which 
contradicts (2) above.
\epf

}

\epF{Theorem~\ref{mt}, \ela\ case}

\punk{Final remarks}
\las{re} 

\bvo
\lam{vo1}
It is interesting to figure out whether 
Theorem~\ref{mt} holds in other extensions of 
$\rL$ by a single generic real, \eg\ in 
extensions by a single Cohen-generic\snos
{Since adding a single Cohen reals is equivalent 
to adding many Cohen reals, it is fairly easy to show 
that there are indiscernible sets of reals 
in Cohen extensions, \eg\ $[a]_\rL$ and $[b]_\rL$ 
for any Cohen-generic pair of reals $\ang{a,b}$, 
as shown in Theorem 3.1 of \cite{ena}.  
On the other hand, such indiscernibles 
hardly form an \OD\ pair, or, equivalently, 
arise as equivalence classes of an OD \eqr\  
$\rE$ with only two equivalence classes.},
or a single Solovay-random, or a single Silver real.
The random case is espesially interesting as 
it is close to the Sacks case in some forcing 
details like the property of \kb\  reading 
of names of reals.
One of the technical difficulties is to prove an 
analog of \gol{\cor ing lemmas in section} \ref{corS} 
for perfect sets of positive measure. 
The merger of equivalence classes, rather 
transparent in the proof of 
lemmas \ref{s2}, \ref{s3}, 
becomes way more complex then. 
\gol{
On the positive side, it turns out that Theorem~\ref{mt} 
also holds for forcing by perfect non-$\sg$-compact sets 
in $\bn,$ to be published elsewhere. 
}
\evo

\bvo
\lam{vo2}
In view of Theorem~\ref{mt}, 
one may ask whether there is a model in which 
every finite non-empty \OD\ set contains an 
\OD\ element but there are
non-empty (infinite) OD sets containing no OD elements. 
Could the Solovay model \cite{sol} 
(where all projective sets are measurable) 
be such a model?
\evo

\punk{History of this result}
\las{his} 

The proof of Theorem~\ref{mt} given above 
was manufactured by V.\,Kanovei in 
January 2020, after a short discussion at 
Mathoverflow\snos
{%
\url{https://mathoverflow.net/questions/349243} }, 
%
on the 
basis of the following exerpt from an email message 
from R.\,M.\;Solovay to Ali Enayat, 
quoted here thanks to 
Solovay's generous permission. 
\begin{quote}
[\ubf Solovay to Enayat 25.10.2002:\rm] 
{

Here's a freshly minted theorem.

Consider the Sacks extension of a model of 
$\rV=\rL +\ZFC$. 
Then LA does not hold.\snos
{In the context of this exchange, LA is the 
Mycielski axiom, the axiom formulated by Mycielski, 
investigated in Enayat's paper \cite{ena}, 
in which it is referred to as the Leibniz-Mycielski 
axiom LM. 
LM states that given any pair of distinct sets 
$a$ and $b$, there is some ordinal $\alpha$, 
and some first order formula $\phi(x)$, 
such that $\rV_\alpha$ contains $a$ and $b$, 
and $\rV_\alpha$ 
satisfies  $\phi(a)$ but does not satisfy  
$\phi(b)$. 
The motivation for establishing Theorem 1.1  
was the guess 
(privately communicated by Enayat to Solovay) 
that the consistency of 
$\ZFC + \text{LM}  +  
\text{``\msur$\rV\ne \text{HOD}$''}$ 
can be shown 
by verifying that LM holds in the 
extension of the constructible universe by a Sacks real.  
The question of consistency of  
$\ZFC + \text{LM}  +  
\text{``\msur$\rV\ne \text{HOD}$''}$ 
has proved 
to be more difficult than meets the eye, 
and remains open.}

My proof is a bit involved. 
Here's a high level - view. 

     By a transfinite construction of length 
$\aleph_1$ I
construct a $P$-name $E$ such that the following are
forced:
\bit
\item
$E$ is an equivalence relation on the set of
non-constructible reals.
\item
$E$ has precisely two equivalence classes.

\item     
In each perfect set with constructible code there
are representatives of both equivalence classes.

\item
$E$ is ordinal definable.
\eit
The two distinct but indiscernable members of
the generic extension are the two equivalence classes
of $E$.
     
The proof is a bit too involved to type in using
a web-interface like yahoo. 
(Shades of Fermat's margin!) 
The proof uses one standard but relatively
deep fact from descriptive set theory. 
If $B$ is an
uncountable Borel set, then $B$ contains a perfect
subset.
     
-- Bob

P.S. 
I don't use much about $\rL$. 
Just that it satisfies
$\rV = \OD$ and is uniformly definable 
in any extension and
that it satisfies CH.\snos
{This is equally true for our proof.} 
}
\ [\ubf End\rm]
\end{quote}
The above proof of Theorem~\ref{mt}  
in the Sacks case 
obviously more or less follows 
Solovay's outline. 
In light of the key role of the Silver Dichotomy 
in the proof presented here, 
we don't know to what degree 
it coincides with the original proof 
by Solovay in all important details.  

Upon the completion of the proof, the co-authors 
contacted R.\,M.\;Solovay, with an invitation to 
join as a co-author of this note, 
but he unfortunately did not accept our invitation.

%

\bibliographystyle{plain}
\addcontentsline{toc}{subsection}{\hspace*{5.5ex}References}
{\small
\bibliography{54}
}

\end{document}